\documentclass[12pt,a4paper]{article}

\title{The Gelfand map and symmetric products }
\author{V.M.Buchstaber and E.G.Rees}
\date{}

\begin{document}

\maketitle

There are many instances of the principle that if $A$ is an algebra of
functions on $X,$ then every ring homomorphism $A \rightarrow {\bf C}$
is given by evaluation at a particular $x \in X.$ Examples are the
nullstellensatz in algebraic geometry and the result of Gelfand when
$X$ is a compact Hausdorff space. In these cases one can regard $X$ as
being included in $ {\rm Hom}(A, {\bf C})$ as the set of those $f : A
\rightarrow {\bf C}$ which satisfy the set of equations $f(xy) =
f(x)f(y)$ indexed by $(x,y) \in X \times X.$ In this paper we
introduce the corresponding equations for the symmetric products of
$X.$ We show that, in these examples, ${\rm Sym}^n(X)$ is included in
${\rm Hom}(A, {\bf C})$ as the set of those $f$ that satisfy these
more complicated equations.  

Given a linear map $f : A \rightarrow {\bf C}$ we consider certain
maps which can be regarded as ``higher'' versions of $f$ and are
denoted $\Phi_n(f) : A^{\otimes n} \rightarrow {\bf C},$ their
definition is based on formulae used by G. Frobenius \cite{[Fro1]}.
The subset $\Phi_n(A) \subset \textup{Hom}(A,{\bf C}) $ of the space
of all linear maps $\textup{Hom}(A,{\bf C})$ consisting of those $f$
for which $\Phi_{n+1}(f)=0 $ and $f(1)=n$ is particularly interesting
and we will develop its properties; $\Phi_1(A)$ is the set of ring
homomorphisms. When $A$ is an algebra of functions on a space $X$ the
sets $\Phi_n(A)$ are closely related to the symmetric product ${\rm
Sym}^n(X):= X^n/\Sigma_n.$ The case $n=1$ is classical, if $A$ is a
separating algebra of functions on a compact space $X$ then
$\Phi_1(A)$ is the set of algebra homomorphisms and, by the Gelfand
transform, this is homeomorphic to $X.$

In many cases, we can identify $\Phi_n(A)$ with the set of maps that
can be written as the sum of $n$ ring homomorphisms. In particular,
when $X$ a compact Hausdorff space and $A=C(X)$ is the ring of complex
valued continuous functions on $X$, $\Phi_n(A)$ is precisely the set
of linear maps that can be written as the sum of $n$ ring
homomorphisms and so can be identified with ${\rm Sym}^n(X)$.  The
analogous result holds when $A$ is a finitely generated commutative
algebra. In particular, when $X = {\bf C}^m$ and the algebra is the
ring $A={\bf C}[u_1,u_2, \dots,u_m]$ of polynomial functions on $X$ we
prove that $\Phi_n(A)$ is the symmetric product ${\rm Sym}^n({\bf
C}^m).$ A by product is that the embedding $\Phi_n(A) \subset
\textup{Hom}(A,{\bf C}) $ can be described by specific equations, we
study them in another paper.

 The equations that define $\Phi_n$ can be described in several ways
and can be derived from formulae introduced by Frobenius \cite{[Fro1]}
to define the $k$-characters of a finite group.  These formulae have
also been used more recently (eg \cite{[For]} ) in the study of
relations in matrix algebras. In the context of these works, the case
where $A$ is non-commutative was of primary interest; indeed such
functions were often trivial in the commutative case.  In contrast, in
this paper we concentrate our attention on the commutative case.

Our initial interest in these ideas came from \cite{[BR1]}
\cite{[BR2]} where the diagonal of an $n$-Hopf algebra (the analogue
of a Hopf algebra for an $n$-valued group) is characterised as being
an $n$-ring homomorphism.  Professor John McKay kindly pointed out to
us that there were similarities to the formulae introduced by
Frobenius, and indeed there is a precise relationship \cite{[BR3]}.

\bigskip
The paper is divided into sections:

\begin{enumerate}
\item Proves an identity satisfied by partitions of sets and which may
be of independent interest.

\item Introduces Frobenius transformations and develops some of their
basic properties. 

\item Sets up the relationship with symmetric products and proves the
basic theorem for finite sets, affine varieties and  compact
Hausdorff spaces.

\end{enumerate}

\bigskip
\begin{center}
{\bf \S 1 An identity on partitions }
\end{center}

If $\sigma$ is a permutation of a set $X$, there is a partition of $X$
given by the orbits of the action of the group generated by $\sigma.$
Clearly, two permutations giving the same partition have the same
cycle type and hence the same sign; given a partition $\pi,$ let
$\epsilon (\pi)$ denote this sign and let $n(\pi)$ be the number of
permutations that give rise to $\pi,$ so if the parts of $\pi$ are
$P_1,P_2,\dots, P_k$ then $n(\pi) = \prod_{i=1}^k(\#P_i-1)!.$

Let ${\cal P}(X)$ denote the free abelian group on the set of all
partitions of $X$ and $$\chi(X) = \sum \epsilon (\pi )n( \pi ) \pi \in
{\cal P}(X)$$ where the sum is over all the partitions of $X.$

If $ \pi_1, \pi_2$ are partitions of $X,Y$ respectively, then one has
a natural partition $ \pi_1 \pi_2$ of the disjoint union $ X \sqcup
Y.$ So, one can define $\chi(X)\chi(Y) \in {\cal P}(X\sqcup Y).$

If $g : X \rightarrow Y$ is a map and $\pi$ is a partition of $Y$
one has, by taking inverse images of the parts of $\pi$, a partition
$g^{\ast}\pi$ (they have the same number of parts if $g$ is
surjective) and hence a homomorphism $g^{\ast} : {\cal P}(Y)
\rightarrow {\cal P}(X).$

\bigskip \noindent {\bf (1.1) Definition }

A {\bf partial pairing} $\phi$ between two sets $X,Y$ is a bijection
$\phi : X_{\phi} \rightarrow Y_{\phi}$ between subsets $X_{\phi}
\subset X$ and $Y_{\phi } \subset Y.$ Given a partial pairing $\phi$,
define an equivalence relation on $X \sqcup Y$ by $x \sim y$ if
$\phi(x)=y.$ We denote the quotient set by $X \sqcup_{\phi} Y$ (its
cardinality is $\# X+ \# Y-\# X_{\phi}$) and the quotient map by
$q_{\phi} : X \sqcup Y \rightarrow X \sqcup_{\phi} Y.$

\bigskip \noindent {\bf (1.2) Proposition }

{\it If $X,Y$ are disjoint, then 

$$\sum q_{\phi}^{\ast}\chi( X \sqcup_{\phi} Y) = \chi(X)\chi(Y)$$
where the sum is over all partial pairings $\phi$ between $X$ and $Y,$
including the `empty' partial pairing. }

\medskip \noindent \underline{Proof} 

\noindent

Let $\pi = P_1 \sqcup P_2 \sqcup \dots  \sqcup P_k$  be a
partition of $X \sqcup Y$ where 

\noindent
$P_j = \{x_{1j},
x_{2j}, \dots ,x_{m_jj}, y_{1j}, y_{2j}, \dots ,y_{n_jj}\}$ with
$x_{ij} \in X, y_{ij} \in Y$ and let
$$c(m,n,\ell)=(-1)^{m+n-\ell-1}\;(m+n-\ell-1)!\;\ell!\;\;{m\choose
\ell}{n \choose \ell}.$$ Then the coefficient of the partition $\pi$
that appears from terms arising from pairing along subsets of cardinality
$\ell$ is
$$\sum_{\ell_1+\ell_2+\dots+\ell_k=\ell} c(m_1,n_1,\ell_1)
c(m_2,n_2,\ell_2) \dots c(m_k,n_k,\ell_k). $$ Hence, the coefficient
of $\pi$ in the left hand side of the expression in Proposition (1.2) is

$$c_{\pi} =\sum_{\ell_1=0}^{\textup{\footnotesize min}(m_1,n_1)}
\sum_{\ell_2=0}^{\textup{\footnotesize min}(m_2,n_2)} \dots
\sum_{\ell_k=0}^{\textup{\footnotesize min}(m_k,n_k)}c(m_1,n_1,\ell_1) c(m_2,n_2
,\ell_2)
\dots c(m_k,n_k,\ell_k). $$

To evaluate this sum we let $$d(m,n,\ell)= \frac{{m \choose \ell}{n \choose
\ell}}{{m+n-1 \choose \ell}} =
(-1)^{m+n-1} (-1)^{\ell}\frac{c(m,n,\ell)}{(m+n-1)!}$$  and
$$P_{m,n}(t)=
\sum_{\ell=0}^{\textup{\footnotesize min}(m,n)} d(m,n,\ell)(-t)^\ell \;;$$
then  $$c_{\pi}= \epsilon(\pi) \prod_{r=1}^{r=k}(m_r+n_r-1)! P_{m_r,n_r}(1).$$
The polynomial  $P_{m,n}(t)$
is a hyper-geometric polynomial, being a solution
of the differential equation
$$t(1-t)y^{''}(t)-(m+n-1)(1-t)y^{'}(t)-mny(t)=0.$$  Therefore, by
substituting $t=1$ into the differential equation
 one has $P_{m,n}(1)=0$ for $\textup{ min}(m,n) >0.$

Hence, $c_{\pi} =0$ unless each part of $\pi$ consists either entirely
of $x$'s or entirely of $y$'s. In this case, let the partition be
$\pi_1\pi_2$; its coefficient in the term $\chi(X)\chi(Y)$ is
$\epsilon(\pi_1)n(\pi_1)\epsilon(\pi_2)n(\pi_2)$ (it can only occur
once in the product). In the left hand side of the equation in
Proposition (1.2) the only term in which $\pi_1\pi_2$ appears is $\chi(X
\sqcup Y)$ and its coefficient is $\epsilon(\pi_1\pi_2)n(\pi_1\pi_2) =
\epsilon(\pi_1)n(\pi_1)\epsilon(\pi_2)n(\pi_2).$

The result of Proposition (1.2) now follows.

\begin{center}
{\bf \S 2 Frobenius transformations }
\end{center}

In \cite{[Fro1]} and \cite{[Fro2]}, G. Frobenius introduced the
$k$-characters of a finite group and they have been studied again more
recently (e.g. \cite{[Joh]}, \cite{[HJ]}). In \cite{[BR3]} we extended
the concept to a broader context and called them Frobenius
transformations. Here we consider the case where the algebras on which
they are defined and in which they take values are both
commutative. These transformations have interesting applications in
these case despite the fact that the $k$-characters considered by
Frobenius vanish for irreducible representations of finite abelian
groups.

Following \cite{[For]}, where the Frobenius formula for the
$k$-character is reformulated, write a permutation $\sigma \in
\Sigma_{n+1}$ as a product of disjoint cycles (including those of
length one) $$\sigma = \gamma_1 \; \gamma_2 ...\;\gamma_q.$$ If $f : A
\rightarrow B$ is linear and $\gamma$ is the cycle
$(r_1~r_2~\dots~r_k)$ write $f_{\gamma}(a_1,a_2,...a_{n+1})
=f(a_{r_1}a_{r_2}...a_{r_k}).$ (Note that this also works for
non-commutative algebras $A$ provided $f$ is `trace-like' because in
that case the value of $f_{\gamma}(a_1,a_2,...a_{n+1})$ depends only
on the cycle $\gamma$ and is independent of the way that it is written
in terms of the $x$'s.) Now write
$$f_{\sigma}=f_{\gamma_1}\;f_{\gamma_2}\;...\;f_{\gamma_q}\;.$$

\bigskip \noindent {\bf (2.1) Definition } 
{\it For a linear map $f : A \rightarrow B$ where $A,B$ are commutative
algebras, the map $\Phi_m(f) : A^{\otimes m}\rightarrow B$ is defined
by }
$$\Phi_m(f)(a_1,a_2,\dots, a_m) = \sum_{\sigma \in
\Sigma_m}\epsilon(\sigma)f_{\sigma}(a_1,a_2,\dots, a_m).$$

\bigskip \noindent {\bf Remark } The map  $\Phi_m(f)$ is clearly
symmetric and multilinear.

\bigskip We can define $ f(\chi (X))$ as follows :

\bigskip \noindent {\bf (2.2) Lemma } 
{\it If $X = (a_1,a_2,\dots, a_m)$ and $\chi(X)$ is as defined in \S 1
then $\Phi_m(f)(a_1\otimes a_2\otimes \dots a_m) = f(\chi (X)).$}

\medskip \noindent \underline{Proof} 

\noindent

Rewriting the above definition, if $P =
[a_{i_1},a_{i_2},\dots,a_{i_r}]$ is a multi-subset of the multi-set $ X
=[a_1,a_2,\dots, a_m ] \subset A $ then $f_P(a_1,a_2,\dots, a_m) =
f(a_{i_1}a_{i_2}\dots a_{i_r})$ and if $\pi = P_1 \sqcup P_2 \sqcup
\dots \sqcup P_k$ is a partition of $X$, we have $f_{\pi } =
f_{P_1}f_{P_2}\dots f_{P_k}$ and $\Phi_m(f) = \sum \epsilon(\pi
)n(\pi)f_{\pi }.$

\medskip
We now recall the inductive definition based on that used by Frobenius. 

\bigskip \noindent {\bf (2.3) Definition } 
{\it 
 Define, inductively, for $n \in {\bf N}$, linear maps $\Phi_n(f) :
A^{\otimes n} \rightarrow B$ starting with $\Phi_1(f) =f,
\;\;\Phi_2(f)(a_1,a_2)= f(a_1)f(a_2)-f(a_1a_2)$ and for $n \geq 2$ as
follows :}
$$ \Phi_{n+1}(f)(a_1,a_2,\dots,a_{n+1}) =
f(a_1)\Phi_n(f)(a_2,a_3,\dots,a_{n+1}) - \hspace{2in}$$
$$\Phi_n(f)(a_1a_2,\dots,a_{n+1})-
\Phi_n(f)(a_2,a_1a_3,\dots,a_{n+1})-\dots -
\Phi_n(f)(a_2,a_3,\dots,a_1a_{n+1 }).$$

\medskip \noindent
It is a simple consequence of Proposition (1.2) with $X=\{a_1\},
Y=\{a_2,a_3,\dots , a_{n+1}\}$ that the two definitions (2.1) and
(2.3) are the same.

\medskip \noindent {\bf Remark} It follows immediately from the
inductive definition that if $f$ satisfies  $\Phi_{n}(f) \equiv
0$ then $\Phi_{n+1}(f) \equiv 0.$

\medskip \noindent
{\bf (2.4) Lemma } 
{\it If $B$ is a domain and 
$ \Phi_{n+1}(f) \equiv 0$ but $\Phi_n(f) \not\equiv 0$
then $f(1)=n.$}

\medskip \noindent
\underline{Proof} 

\noindent 
Let $a_1 =1$ then, using the inductive definition, we get 
$$0=\Phi_{n+1}(f)(1,a_2,a_3,...,a_{n+1})= 
[f(1)-n]\Phi_n(f)(a_2,a_3,...,a_{n+1}).$$

\noindent
But, since $\Phi_n(f) \not\equiv 0$ , there are $a_2,a_3,...,a_{n+1}
\in A$ such that 

\noindent
$ \Phi_n(f)(a_2,a_3,...,a_{n+1}) \neq 0.$

\medskip \noindent
{\bf (2.5) Corollary} 
{\it If $ f : A \rightarrow B$ satisfies   $\Phi_{n+1}(f) \equiv 0$  and 
$B$ is a domain then $f(1) \in \{0,1,2,...,n\}.$ }

\medskip \noindent
\underline{Proof}  

\noindent This follows by a simple induction from Lemma (2.4).

\medskip \noindent
{\bf (2.6) Definition} 
{ \it A linear map $f : A \rightarrow B$ is a  
{\bf Frobenius $n$-homomorphism } if $\Phi_{n+1}(f) \equiv 0$
and $f(1)=n.$}

\medskip \noindent
{\bf Remark }
When $B$ is a domain, every Frobenius 1-homomorphism 
$ f : A \rightarrow B$ is a ring homomorphism.

\medskip \noindent {\bf (2.7) Proposition}{ \it If $B$ is a domain
then a linear map $f : A \rightarrow B$ such that $\Phi_{n+1}(f)
\equiv 0$ and  $f(1)=k \leq n$ is a Frobenius $k$-homomorphism .}

\medskip \noindent
\underline{Proof}

 \noindent Applying the inductive formula to
 $\Phi_{n+1}(f)(1,a_2,\dots,a_{n+1}) =0$ gives,

$$(k-n)\Phi_{n}(f)(a_2,\dots,a_{n+1})=0$$ 
 \noindent
so $f$ is a Frobenius  $k$-homomorphism. The result follows by induction.

\bigskip 

We denote the sub-algebra of symmetric tensors in $ A^{\otimes n}$ by
${\mathcal S}^nA.$ The map $ \Phi_{n}(f)/n!$ restricted to ${\mathcal
S}^nA$ has a multiplicative property:

\bigskip\noindent
{\bf (2.8) Theorem}
{\it If $f : A \rightarrow B$ is a Frobenius $n$-homomorphism, then the map
defined by $$\frac{\Phi_n(f)}{n!} : {\mathcal S}^nA\rightarrow B$$ is a ring
homomorphism.}

\medskip \noindent
\underline{Proof} 

A typical element of ${\mathcal S}^nA$ is $${\mathbf a}= \sum_{\sigma \in
\Sigma_n}a_{\sigma(1)}\otimes a_{\sigma(2)}\otimes \dots \otimes
a_{\sigma(n)}$$ and so the product of two such elements is 

$$ {\mathbf a}{\mathbf b} =\sum_{\sigma_1,\sigma_2 \in
\Sigma_n}a_{\sigma_1(1)}b_{\sigma_2(1)}\otimes
a_{\sigma_1(2)}b_{\sigma_2(1)}\otimes \dots \otimes
a_{\sigma_1(n)}b_{\sigma_2(n)}.$$

 \noindent
By Lemma (2.2), if $X= (a_1,a_2,\dots, a_n), Y= (b_1,b_2,\dots, b_n)$ we
have $$\Phi_n(f)(a_1, a_2, \dots ,a_n)\Phi_n(f)(b_1,
b_2, \dots ,b_n)= f(\chi(X))f(\chi(Y)).$$

 \noindent
Regarding $X,Y$ as disjoint we have  $$
f(\chi(X))f( \chi(Y)) = f(\chi(X) \chi(Y)).$$
 \noindent
 By Proposition (1.2) one has
$$f(\chi(X) \chi(Y))= \sum_{\phi \in
\Sigma_n}f(q_{\phi}^{\ast}\chi(X\sqcup_{\phi} Y)$$ since $f$ is a
Frobenius $n$-homomorphism the only terms on the right hand side which
are non-zero are those where $\phi : X \rightarrow Y $ is a
permutation.

 \noindent By Lemma (2.2), $$\sum_{\phi \in
\Sigma_n}f(q_{\phi}^{\ast}\chi(X\sqcup_{\phi} Y))= \sum_{\phi \in
\Sigma_n}\Phi_n(f)(a_1b_{\phi(1)}, a_2b_{\phi(2)}, \dots
,a_nb_{\phi(n)}). $$

 \noindent
Hence adding all the relevant terms we get 
$$\Phi_n(f)({\mathbf a})\Phi_n(f)({\mathbf b})  \hspace{9cm}$$
$$=\sum_{\sigma_1,\sigma_2,\sigma_3 \in
\Sigma_n}\Phi_n(f)(a_{\sigma_1(1)}b_{\sigma_2(1)\sigma_3(1)},\;
a_{\sigma_1(2)}b_{\sigma_2(2)\sigma_3(2)}, \dots
,a_{\sigma_1(n)}b_{\sigma_2(n)\sigma_3(n)} ) \hspace{1cm}$$
$$=n!\sum_{\sigma_1,\sigma_2 \in
\Sigma_n}\Phi_n(f)(a_{\sigma_1(1)}b_{\sigma_2(1)},\;
a_{\sigma_1(2)}b_{\sigma_2(2)}, \dots
,a_{\sigma_1(n)}b_{\sigma_2(n)} )=n!\Phi_n(f)({\mathbf a}{\mathbf b})$$

\noindent
It is easy to check that $\Phi_n(f)(1, 1, \dots ,1)=n!.$

\medskip \noindent
{\bf (2.9)  Theorem} {\it If $f,g$ are Frobenius $m$- and $n$-homomorphisms
respectively, then $f+g$ is a Frobenius  ($m+n$)-homomorphism. }

\medskip \noindent
{\bf (2.10)  Corollary}  {\it If $f : A \rightarrow
B$ is the sum of $n$ ring homomorphisms $f_i : A \rightarrow B,\; 1
\leq i \leq n$ then $f$ is a Frobenius $n$-homomorphism.}

\medskip \noindent To prove (2.9) it is convenient to use polarisation
(\cite{[W]}, Chapter II) to study the properties of $ \Phi_{n}(f);$
this means that because $\Phi_n(f)(a_1,a_2, \dots,a_n)$ is both
multi-linear and symmetric it is enough to calculate using `diagonal'
elements only.

\medskip
 For a partition $\lambda =\{\lambda_1,\lambda_2,\dots, \lambda_q\}$
of $n$ we let $$f_{\lambda}(a,a,
\dots,a)=f(a^{|\lambda_1|})f(a^{|\lambda_2|}) \dots
f(a^{|\lambda_q|})$$ and $\epsilon (\lambda)$ is the sign of a
permutation whose cycle decomposition consists of cycles of lengths
$\{\lambda_1,\lambda_2,\dots, \lambda_q\}$. Hence, from (2.1)

$$\Phi_{n}(f)(a,a, \dots,a) =\sum_{\lambda} {\epsilon}(\lambda)
n(\lambda)f_{\lambda}(a,a, \dots,a)$$ where $n(\lambda)$ denotes the number of
elements of the symmetric group $\Sigma_n$ in the conjugacy class
determined by $\lambda.$

\medskip \noindent {\bf (2.11) Lemma}
$$\Phi_{n}(f)(a,a, \dots,a)=(n-1)!\sum_{k=1}^{n}(-1)^{k+1}f(a^k)
\frac{\Phi_{n-k}(f)(a,a, \dots,a)}{(n-k)!}$$

\medskip 
\noindent
\underline{Proof} 

This is obtained from Definition (2.1) by breaking the sum 

$$ \Phi_{n}(f)(a,a, \dots,a) = \sum_{\sigma \in
\Sigma_{n}}\epsilon(\sigma)f_{\sigma}(a,a, \dots,a) $$ into parts
corresponding to the length of the cycle in the permutation $\sigma$
that contains $n$ and there are $(n-1)(n-2) \dots (n-k+1)$ such cycles
of length $k$.

\medskip \noindent
{\bf (2.12)  Corollary} {\it The exponential generating function

$$\sum_{n=0}^{\infty}\frac{\Phi_{n}(f)(a,a, \dots, a)}{n!}t^n =
\exp\left(\sum_{k=1}^{\infty}(-1)^{k+1}\frac{f(a^k)}{k}t^k\right)$$}

\medskip 
\noindent
\underline{Proof} 
This follows from the following well known combinatorial result.

\medskip 

\noindent {\bf (2.13) Lemma} {\it If $\; \; \Phi_0=1$ and $$\Phi_n =
(n-1)!\sum_{k=1}^n\frac{s_k\Phi_{n-k}}{(n-k)!}\;\; {\rm for}\;\; n\geq 1$$
then $$\sum_{n=0}^{\infty}\frac{\Phi_n}{n!}t^n =
\exp\left(\sum_{k=1}^{\infty}\frac{s_k}{k}t^k\right).$$}

\medskip 
\noindent
\underline{Proof} 

\noindent
 $${\rm Let }\;\; \Phi(t)= \sum_{ n=0}^{\infty}\frac{\Phi_n}{n!}t^n\;\;{\rm and} \;\; s(t)= \sum_{k=0}^{\infty}s_kt^k
\hspace{2.5in}$$ 
then
the hypothesis gives that
$$t\Phi^{\prime}(t) = \Phi(t)s(t)$$
so $$\log(\Phi(t)) = \int\frac{s(t)}{t}$$
Checking the constant term gives the required conclusion.

\medskip 
\noindent
\underline{Proof of (2.9)} 

We first note that by definition a map $f$
is a Frobenius $n$-homomorphism if and only if the exponential
generating series is a polynomial of degree $n$ and $f(1)=n.$ Using
(2.12) we get that

$$\sum_{k=0}^{\infty}\frac{\Phi_{k}(f+g)(a,a,\dots ,a)}{k!}
=\sum_{r=0}^{m}\frac{\Phi_{r}(f)(a,a,\dots
,a)}{r!}\sum_{s=0}^{n}\frac{\Phi_{s}(g)(a,a,\dots ,a)}{s!},$$ so
$\Phi_{r}(f+g)(a,a,\dots ,a)=0$ for $r > m+n$ and clearly $(f+g)(1)=m+n.$

\smallskip
We will be interested in the behaviour of Frobenius $n$-homomorphisms
on idempotent elements in $A,$ and the following result is a
generalisation of Corollary (2.5).

  \medskip \noindent
{\bf (2.14) Lemma } 
{\it If $B$ is a domain and $a \in A$ satisfies $a^2=a$ and 

\noindent
$ \Phi_{n+1}(f)(a,a,\dots ,a)=  0,$ 
then $f(a)=k$ for some integer $k$ such that

\noindent
$0 \leq k \leq n.$}

\medskip \noindent
\underline{Proof} 

\noindent
If $n=1,$ then since $ \Phi_2(f)(a,a) =0,$ 
$f(a)^2=f(a^2)$ and so 
$f(a)(f(a)-1)=0.$ 

\noindent
If $ \Phi_{n+1}(f)(a,a,\dots
,a)= 0,$ and $a^2=a$  then  $$f(a) \Phi_n(f)(a,a,\dots ,a) - n
\Phi_n(f)(a,a,\dots ,a) = 0.$$ So, either, $f(a)=n$ or $
\Phi_n(f)(a,a,\dots ,a)=0$ and the result follows by induction. 

 \medskip \noindent
{\bf (2.15) Lemma } 
{\it If $a \in A$ then

\noindent
$ \Phi_{n}(f)(a,1,\dots ,1)=  f(a)(f(1)-1)(f(1)-2) \dots (f(1)-(n-1)).$ }

\medskip \noindent
\underline{Proof} 

This is straightforward to check for $n=1,2$ and then by induction.

\bigskip

\begin{center}

{\bf \S 3 Symmetric products}

\end{center}

In this section we relate our theory to the study of 
symmetric products.

\bigskip
The set of all Frobenius $n$-homomorphisms from $A$ to ${\bf C}$ will
be denoted by $\Phi_n(A).$

Although we will prove a more general version of the following result,
it is worth starting with the following simple  proof.

\bigskip\noindent
{\bf (3.1) Theorem} {\it
For a finite set $X,$ the evaluation map

$${\cal E} : {\rm Sym}^n(X) \rightarrow \mathrm{Hom}(C(X), \mathbf{C})$$

 \noindent
defined by $[x_1, x_2, \dots, x_n] \rightarrow \{f \rightarrow \sum f(x_r)\}$
is an isomorphism onto the set of  Frobenius $n$-homomorphisms.}

\medskip \noindent \underline{Proof} 

\noindent
First note that evaluation at a point is a ring homomorphism $C(X)
\rightarrow \mathbf{C}$ and so ${\cal E}[x_1,x_2,\dots,x_n]$ is a
Frobenius $n$-homomorphism by Corollary (2.10).

\noindent
Let $D=[x_1, x_2, \dots, x_n] \in {\rm Sym}^n(X)$ be considered as a formal sum
$D=\sum m_rx_r,$ where $m_{r} \in \mathbf{Z}_{+}$ and $\sum m_r =n$.
Let $e_r \in C(X)$ be the function that is 1 on $x_r$ and 0 elsewhere.
Then ${\cal E}(D)e_r = m_r.$ 

Now if $f : C(X) \rightarrow \mathbf{C}$ is a Frobenius
$n$-homomorphism and not an $(n-1)$-homomorphism, then, by 
(2.14), $f(e_r) = f_r \in \mathbf{Z}_{+}$ and because $1$ is the sum
of the idempotents $e_r, \sum f_r =n.$
Then $D= \sum f_rx_r$ is mapped onto $f$ by ${\cal E}.$

The map ${\cal E}$ is injective for a very general class of algebras
of functions on a space $X$
(those where functions separate points) and so it is an isomorphism in
the case of finite sets.

\bigskip \noindent {\bf (3.2) Theorem } 

{\it If $f : {\bf C}[u_1,u_2,\dots, u_m] \rightarrow {\bf C}$ is a
Frobenius $n$-homomorphism then there are points ${\bf x}_1,{\bf
x}_2,\dots {\bf x}_n \in {\bf C}^m$ such that $f(p)=p({\bf x}_1)+p({\bf
x}_2)+\dots+ p( {\bf x}_n).$}

\medskip \noindent \underline{Proof} 

\noindent

By Theorem (2.8), the map $$\frac{\Phi_n(f)}{n!} : {\mathcal S}^n({\bf
C}[u_1,u_2,\dots, u_m]) \rightarrow {\bf C}$$ is a ring
homomorphism. But ${\mathcal S}^n({\bf C}[u_1,u_2,\dots, u_m])$ is the
algebra of polynomial functions on $\textup{Sym}^n({\bf C}^m)$, so
$\frac{\Phi_n(f)}{n!}$ is given by evaluation at a multi-set $[{\bf x}_1,{\bf
x}_2,\dots, {\bf x}_n] \subset {\bf C}^m.$ If $p \in {\bf
C}[u_1,u_2,\dots, u_m]$ then, on the one hand, by Theorem (2.8),
 $$\frac{\Phi_n(f)}{n!} ((p,1,1, \dots,1)+(1,p,1, \dots,1) + \dots
+(1,1,\dots,p)) =p(x_1)+p(x_2)+ \dots +p(x_n)$$ and, on the other hand,
by the definition of $\Phi_n(f)$, the fact that $f(1)=n$ and Lemma
(2.15) one gets that $$\frac{\Phi_n(f)}{n!} ((p,1,1, \dots,1)+(1,p,1, \dots,1)
+ \dots +(1,1,\dots,p)) = f(p).$$

This proves Theorem (3.2).

\bigskip

Reformulating this result and denoting the set of all Frobenius
$n$-homomorphisms $f : {\bf C}[u_1,u_2,\dots, u_m] \rightarrow {\bf
C}$ by $\Phi_n({\bf C}^m)$ gives

\bigskip 
\noindent {\bf (3.3) Corollary } 

{\it The evaluation map $${\cal E} : \textup{Sym}^n({\bf C}^m) \rightarrow
\Phi_n({\bf C}^m)$$ is a homeomorphism.}

\bigskip 

\noindent {\bf (3.4) Theorem } 

{\it Let $A$ be a finitely generated commutative algebra and let 
$f : A  \rightarrow {\bf C}$ be a Frobenius $n$-homomorphism then
there are ring homomorphisms $f_i : A  \rightarrow {\bf C}$ for 
$ 1 \leq i \leq n$ such that $f = f_1 + f_2 + \dots + f_n.$}

\medskip \noindent \underline{Proof} 

\noindent

We deduce this from Theorem (3.2) (which is the special case when $A$ is a
polynomial algebra). Let $q :{\bf C}[u_1,u_2,\dots, u_m] \rightarrow
A$ be a quotient map onto $A$ whose kernel we denote by $I.$ Then $g =
fq$ is a Frobenius $n$-homomorphism on ${\bf C}[u_1,u_2,\dots, u_m]$
and so, by Theorem (3.2), there are ring homomorphisms $g_i : {\bf
C}[u_1,u_2,\dots, u_m] \rightarrow {\bf C}$ with $g= g_1 + g_2 + \dots
+ g_n.$ We show that each $g_i$ vanishes on $I.$

First we relabel the $g_i$ so that $g_1,g_2, \dots, g_k$ are distinct
and $$g= r_1g_1+r_2g_2+ \dots r_kg_k \;\;\textup{for} \;\; r_i \in {\bf N}.$$ If
$\theta$ is a polynomial in $u_1,u_2,\dots, u_m$, then 
$g_i(\theta) = \theta(g_i(u_1),g_i(u_2),\dots, g_i(u_m))$ since $g_i$
is a ring homomorphism, and if $\psi$ is another such polynomial
$g_i(\psi \theta)=g_i(\psi)g_i(\theta).$ If $\theta \in I$  
then $g(\theta)=0$ and $g(\psi \theta)=0$ so 
$$
   \left(
\begin{array}{cccc}
1 & 1 &   \dots &1 \\
g_1(\psi) &g_2(\psi)  &    \dots &g_k(\psi) \\
\vdots &\ddots & \ddots&  \vdots \\
\vdots &\ddots & \ddots&  \vdots \\
g_1(\psi)^{k-1} & g_2(\psi)^{k-1} &  \dots&g_k(\psi)^{k-1}
\end{array} \right) 
\left(
\begin{array}{c}
 r_1g_1(\theta) \\
  r_2g_2(\theta) \\
\vdots  \\
\vdots  \\
  r_kg_k(\theta) 
\end{array} \right) =\left(
\begin{array}{c}
 0 \\
 0 \\
\vdots  \\
\vdots  \\
 0 
\end{array} \right) 
                   $$


\bigskip \noindent {\bf (3.5) Lemma }

{\it Given distinct linear maps    $$g_1,g_2, \dots, g_k:
{\bf C}[u_1,u_2,\dots, u_m] \rightarrow {\bf C},$$
 there is a $\psi \in {\bf C}[u_1,u_2,\dots, u_m]$ such that 
$g_1(\psi),g_2(\psi), \dots, g_k(\psi) \in {\bf C}$ are distinct.}

\medskip \noindent \underline{Proof} 

\noindent

By the assumption that the maps are distinct, Ker$(g_i-g_j)$ has codimension 
one for every pair $i \neq j$. Any $\psi \not\in \bigcup_{i \neq
j}\textup{Ker}(g_i-g_j)$ will do.

\bigskip

For such a $\psi$, the matrix above is non-singular and so 
$r_1g_1(\theta) = r_2g_2(\theta)= \dots = r_kg_k(\theta)=0.$ Hence each
$g_i$ vanishes on $I$ and so defines a map $f_i$ on $A$
such that $f=r_1f_1+r_2f_2+\dots r_kf_k,$ that is $f$ is the sum of
$n$ ring homomorphisms.

\bigskip \noindent {\bf Remark} Hence a set of distinct ring
homomorphisms $\{g_1,...,g_k\}$ gives a set of distinct vectors
$\{v_1,...,v_k\}$ and there is a polynomial $\Psi,$ such that the
numbers $\Psi(v_i)=g_i(\Psi)$ are distinct.

\bigskip \noindent {\bf (3.6) Corollary }

{\it Let $A$ be a finitely generated commutative algebra and $V =
 \Phi_1(A).$ Then the evaluation map $ \textup{ Sym}^n(V) \rightarrow
 \Phi_n(A)$ is an isomorphism of varieties.}

\smallskip
If $X$ is a compact Hausdorff space, let ${\rm Sym}^n(X)$ denote the
symmetric product $X^n/\Sigma_n$ and ${\bf C}(X)$ the algebra of
continuous functions on $X$ then the `evaluation' map
$${\cal E}: {\rm Sym}^n(X) \rightarrow \Phi_n({\bf C}(X),{\bf C})$$
defined by $$[x_1,x_2,...,x_n] \rightarrow \left(g \rightarrow 
\sum g(x_k)\right)$$

\noindent
is an embedding. It is clear that ${\cal E}$ is natural and so, if $X$
admits a group action, ${\cal E}$ is equivariant. Then, using the
Gelfand transform that $ \Phi_1 ({\bf C}(X),{\bf C}) \cong X$ one sees
that this is a special case of the general question whether $ {\rm
Sym}^n(\Phi (A,B)) \rightarrow \Phi_n(A,B)$ is an isomorphism when
$B$ is any commutative domain.

We have already shown that when $X$ is a finite set the map ${\cal E}$
is onto.

\bigskip
\noindent
{\bf (3.7) Theorem } 
{\it If $X$ is a compact Hausdorff space and the
function space
${\bf C}(X)$ has the supremum norm, then the map 
$${\cal E}: {\rm Sym}^n(X) \rightarrow \Phi_n^c({\bf C}(X),{\bf C})$$
is a homeomorphism when
 the space of continuous linear functionals on ${\bf C}(X)$ has the
weak topology.}

\bigskip
\noindent
{\bf (3.8) Corollary  } 
{\it Under these conditions every continuous Frobenius $n$-homomorphism is
the sum of $n$ continuous ring homomorphisms.}

\medskip \noindent \underline{Proof of (3.7)} 

This is an easy adaptation of the proof of (3.4). If $f : C(X)
\rightarrow {\bf C}$ is a continuous Frobenius $n$-homomorphism, then it
is easily checked that $$\Phi_n(f)/n!: {\mathcal S}^n(C(X))\rightarrow
{\bf C}$$ is a continuous ring homomorphism. But $ {\mathcal
S}^n(C(X))$ is isomorphic to the algebra $C({\rm Sym}^n(X))$ and the
result follows as in (3.4).

\bigskip
\noindent {\bf Remark}
The case $n=1$ is the classical Gelfand transform map. However,
it seems that the `standard' proofs that the Gelfand map is an
isomorphism do not adapt to prove this more general result. Indeed
most of the proofs for the classical case $n=1$ do not
`find' the point of $X$ at which the ring homomorphism ${\bf C}(X)
\rightarrow {\bf C}$ is evaluation. In another paper we will give such
a constructive proof of (3.7).

\bigskip

{\bf  Acknowledgements}
\medskip

We are very grateful to John McKay for pointing out the relevant
papers by Frobenius and to Mark Haiman who explained to us that
$\Phi_n$ is the expression for a multinomial symmetric function in
terms of the power sums that we use, which lead to these proofs of the
main theorems. We are also grateful to John Byatt-Smith for a
helpful conversation and to Kostya Feldman for helpful comments on an
earlier draft.

The research on which this is based was mainly carried out during
visits by VMB to the University of Edinburgh and supported by  the
Engineering and Physical Sciences Research Council and the London
Mathematical Society. Some of it was carried out at the Programme on
Singularity Theory at the Isaac Newton Institute, Cambridge.

\pagebreak

\vspace{1cm} 

\begin{tabular}{ll}

 Dept. of Mathematics and Mechanics, &
Dept. of Mathematics and Statistics,\\

 Moscow State University, &        James Clerk Maxwell Building,\\

119899,       & King's Buildings,\\

Moscow,   &      Edinburgh EH9 3JZ, \\

Russia. &        Scotland.\\

       & \\

buchstab@mech.math.msu.su & elmer@maths.ed.ac.uk

\end{tabular}

\end{document}